\documentclass[11pt, leqno]{article}
\usepackage{a4, indentfirst, latexsym,amssymb}
\usepackage{t1enc}
\usepackage{color}

  \newtheorem{Th}{Theorem}[section]
  
  \newtheorem{Prop}[Th]{Proposition}

  \newtheorem{Lem}[Th]{Lemma}

  \newtheorem{Rem}[Th]{Remark}

  \def\d{\,\mathrm{d}}

  \def\ind{\mathrm{ind}}

  \def\Deg{\mathrm{Deg}}
  \def\la{\langle}
  \def\ra{\rangle}
  \def\id{\mathrm{id}}

  \def\o#1{\overline{#1}}

  \def\part{\partial}

  \def\Liminf{\mathop{\mathrm{Liminf}\,}}

  \def\R{\mathbb{R}}

  \def\Gr{\mathrm{Gr}}

  \newcommand{\lma}{\lambda}

  \newcommand{\be}{\begin{equation}}
  \newcommand{\ee}{\end{equation}}

\begin{document}

\begin{center}

{\huge  Positive stationary solutions for $p$-Laplacian problems with nonpositive perturbation}

Aleksander \'{C}wiszewski \footnote{Corresponding author}, \ Mateusz Maciejewski ,
\footnote{ The research supported by the MNiSzW Grant no. N N201 395137.}\\

{\small {\em Nicolaus Copernicus University\\
Faculty of Mathematics and Computer Science\\
ul. Chopina 12/18, 87-100 Toru\'n, Poland }\\
e-mail: Aleksander.Cwiszewski@mat.umk.pl, Mateusz.Maciejewski@mat.umk.pl\\

Version: 05.04.2012.
}
\end{center}

\begin{abstract}
The paper is devoted to the existence of positive solutions of nonlinear elliptic equations with $p$-Laplacian. 
We provide a general topological degree that detects solutions of the problem
$$
\left\{
\begin{array}{l}
A(u)=F(u)\\
u\in M
\end{array} \right. 
$$
where $A:X\supset D(A)\to X^*$ is a maximal monotone operator in a Banach space $X$ and $F:M\to X^*$
is a continuous mapping defined on a closed convex cone $M\subset X$.
Next, we apply this general framework to a class of partial differential equations with $p$-Laplacian under Dirichlet boundary conditions. In the paper we employ general ideas from \cite{Cw-Kr-JDE2009}, where a setting suitable 
for the one dimensional $p$-Laplacian was introduced. 
\end{abstract}

\thispagestyle{empty}

\section{Introduction}

We shall be concerned with solutions to the following nonlinear boundary value problem 
\be\label{29102011-1511} 
\left\{
\begin{array}{l}
-\mathrm{div} (|\nabla u(x)|^{p-2} \nabla u(x)) = f(x,u (x)), \   x\in \Omega,\\
\ \   u(x)\geq 0,\  x\in \Omega,\\
\ \   u(x) = 0,\  x\in \partial \Omega
\end{array} \right. 
\ee
where $\Omega \subset \R^N$ ($N\geq 1$) is a bounded domain with the smooth boundary $\partial \Omega$,
$p\geq 2$ and $f:\Omega\times [0,+\infty) \to \R$ is a Carath\'{e}odory function 
(\footnote{Recall that we say that $f:\Omega\times [0,+\infty)\to\R$ is a Carath\'{e}odory function if 
$f(\cdot, s)$ is measurable for all $s\in [0,+\infty)$ and $f(x,\cdot)$ is continuous for almost all $x\in\Omega$.}).
The differential term  $\mathrm{div} (|\nabla u (x)|^{p-2}\nabla u(x))$ is referred to as the $p$-Laplacian of $u$ at a point $x\in\Omega$. We search for {\em weak solutions} in the Sobolev space $W_{0}^{1,p}(\Omega)$, i.e. $u\in W_{0}^{1,p} (\Omega)$ such that
$$
\int_{\Omega} |\nabla u(x)|^{p-2} \nabla u(x) \cdot \nabla v(x) \d x = \int_{\Omega} f(x, u(x)) v(x) \d x  \ \ \mbox{ for all } \
v\in W_{0}^{1,p} (\Omega).
$$
Such boundary problems with $p$-Laplace were widely studied by many authors who used various methods.
Let us mention just a few. Equations with the one dimensional $p$-Laplacian, i.e. when $N=1$, were studied by Man\'{a}sevich, Njoku i Zanolin \cite{Manasevich-Njoku-Zanolin}, Dr\'{a}bek, Garc\'{i}a-Huidobro and Man\'{a}sevich \cite{Drabek-Huidobro-Manasevich} and as well as by Kryszewski and the author \cite{Cw-Kr-JDE2009}. In the general case, i.e. when $N>1$, positive solutions of $p$-Laplace problems have been studied by a number of authors, e.g. Huang \cite{Huang}, Dr\'{a}bek and Pohozaev \cite{Drabek-Pohozaev}, Ca\~{n}ada, Dr\'{a}bek and G\'{a}mez \cite{Canada-Drabek-Gamez}, Filippiakis, Gasi\'{n}ski and Papageorgiou \cite{Filippiakis-Gasinski-Papa} or Montreanu D., Montreanu V. V. and Papageorgiou \cite{Motreanu-Papa}, V\"{a}th \cite{Vath}.\\
\indent Generally speaking, in the above mentioned papers, either $N=1$ or $N$ is arbitrary but the right has side of the equation - the function $f$ is assumed to be non-negative or satisfy some monotonicity assumptions.
This makes possible to apply Krasnosel'skii's fixed point theorem (in general, fixed point index in cones) or 
variational methods. These assumptions on $f$ seem rather restrictive 
and sometimes unnatural, especially, when we take into account physical interpretation of the considered boundary value problem. In this paper, we do not require $f$ to be non-negative or monotone. 
A general tool for detection of nonnegative solutions is provided. It is based on the geometric idea of tangency 
and using fixed point index in cones. We construct a topological degree for perturbations of maximal monotone operators with respect to closed convex cones. Next we prove appropriate index formulae, which together with the homotopy property, allow us to compute the topological degree in specific examples. 
It is noteworthy, that this setting does not require variational structure and can be also used for systems of $p$-Laplace problems. In this paper, we apply the  method to show the following existence criterion
\begin{Th}\label{29092011-1553}
Suppose that a Carath\'{e}odory function $f:\Omega \times [0,+\infty)\to \R$ and 
$\rho_0, \rho_\infty \in L^\infty (\Omega)$ satisfy the following conditions
\begin{eqnarray}
\label{15112011-2220} & & \mbox{there is } C>0  \mbox{ such that }|f(x,s)| \leq C(1+s^{p-1}) \mbox{ for all } s\geq 0 \mbox{ and a.a. } x\in\Omega;\\
\label{15112011-2222} & & \lim_{s\to 0^+} \frac{f(x,s)}{s^{p-1}} = \rho_0 (x) 
\mbox{ and } \lim_{s\to \infty} \frac{f(x,s)}{s^{p-1}} = \rho_\infty (x) \mbox{ uniformly with respect to } x\in\Omega.
\end{eqnarray}
If the principal eigenvalue $\lma_{1,p}$ of the $p$-Laplace operator lies between $\rho_0$ and $\rho_\infty$, i.e. either
$\rho_0(x) < \lma_{1,p} <\rho_\infty(x)$, for a.a. $x\in\Omega$, or $ \rho_\infty (x) < \lma_{1,p} < \rho_0 (x)$, for a.a. $x\in\Omega$, then the problem {\em(\ref{29102011-1511})} admits a nontrivial weak solution $u\in W_{0}^{1,p} (\Omega)$ such that $u(x)\geq 0$ for a.e. $x\in \Omega$.
\end{Th}
Here the principal eigenvalue $\lambda_{1,p}$ is the smallest real number $\lambda$ such that the problem
\be\label{29092011-1552}
\left\{
\begin{array}{l}
-\mathrm{div} (|\nabla u(x)|^{p-2} \nabla u(x)) = \lma |u (x)|^{p-2} u(x),\   x\in \Omega\\
\ \   u(x)\geq 0,\  x\in \Omega\\
\ \   u(x) = 0,\  x\in \partial \Omega
\end{array} \right. 
\ee
admits a nonzero weak solution (see Remark \ref{04122011-1324} for more details).
Theorem \ref{29092011-1553} corresponds directly to the result of \cite{Huang}, obtained by different methods (the sub-supersolution technique  and the existence result for variational inequalities) and under different assumptions corresponding to the inequality $\rho_\infty < \lma_{1,p} <\rho_0$. Our general method allows us to consider also the case $\rho_\infty> \lma_{1,p}> \rho_0$.
 
The paper is organized as follows. In Section 2  we develop a topological degree detecting coincidence points of maximal monotone operators and continuous operators in closed convex cones. This general tool will be useful if we rewrite the problem (\ref{29102011-1511}) in the form
$$
\left\{ \begin{array}{l} A_p u = N_f (u)\\
u\in M_p
\end{array}\right.
$$ 
where $A_p:L^p(\Omega) \supset D(A_p)\to L^p (\Omega)^*$ is the maximal monotone operator determined by the $p$-Laplacian, $N_f:L^p(\Omega)\to L^p (\Omega)^*$ is the Nemytzkii type operator associated with $f$ and
$M_p$ is the closed convex cone of all non-negative elements in the space $L^p (\Omega)$.
Section 3 provides a general setting in which assumptions of Section 2 are verified. 
Next, in Section 4 we show that the problem (\ref{29102011-1511}) falls into the setting
and, using our topological degree together with spectral properties of $p$-Laplacian, we derive topological index formulae. They turn out to be essential in the proof of Theorem \ref{29092011-1553}, 
which is provided at the end of Section 4.

\noindent {\bf Notation}\\
\indent If $X$ is a metric space  and $B \subset X$, then
$\partial B$ and $cl B$ stand for the boundary of $B$ and the closure of $B$, respectively. If $x_0\in X$ and $r>0$, then
$B (x_0,r):=\{x\in M\mid d(x,x_0)<r  \}$.\\
\indent  If $E$ is a normed space, then by $\|\cdot\|$ we denote its norm.
If $E$ is a normed space and $E^*$ its dual space (of all continuous linear functionals), 
then $\la \cdot, \cdot \ra = \la \cdot, \cdot \ra_{E}:E^*\times E \to \R$ denotes the duality operator 
$\la p,u\ra:=p(u)$, $p\in E^*$, $u\in E$.
If $V$ is another normed space then ${\cal L}(V,E)$ stands for the space of
all bounded linear operators with domain $V$ and values in $E$ with the operator norm denoted by $\| \cdot \|_{{\cal L}(V,E)}$ or simply $\| \cdot \|$
if no confusion may appear.\\
\indent For $x\in\R^N$, $N\geq 1$, $|x|$ denotes the Euclidean norm of $x$ and $x\cdot y$ is the Euclidean scalar product of $x,y\in\R^N$.

\section{Constrained topological degree for perturbations of maximal monotone operators } \label{30112011-2213}

In this section we provide a construction of a topological degree detecting solutions of the abstract constrained problem
\be\label{17102011-1420}
\left\{\begin{array}{l}
0 \in -A u +F(u)\\
u\in M
\end{array} 
\right.
\ee
where $A:X\supset D(A)\multimap X^*$ is a densely defined maximal monotone operator, the constraint set $M$ is a subset of $X$ and  $F: \overline U \to X^*$ is a continuous mapping defined on the closure of an open bounded $U\subset M$.
Throughout the whole section we make the following assumptions\\[0.5em]
\noindent $({\cal A}_1)$ \parbox[t]{140mm}{ there is a 
homeomorphism $N:X\to X^*$ such that $N$ is bounded on bounded sets and the mappings $J_\alpha: X^*\to X$, $\alpha>0$, 
$$
J_\alpha (\tau):= u,  \ \mbox { where }  u\in D(A) \mbox{ is the unique element such that } \tau \in (N+\alpha A)(u),
$$
are well defined and continuous;}\\[0.5em] 
\noindent $({\cal A}_2)$ \parbox[t]{140mm}{ the mapping  ${\cal J}:X^* \times (0,+\infty)\ni (\tau, \alpha )\mapsto J_{\alpha} (\tau )\in X$  is bounded on bounded sets and such that ${\cal J}_{\ |\, X^* \times [\alpha_1, \alpha_2]}$ 
is completely continuous if $0<\alpha_1 \leq \alpha_2$;}\\[0.5em] 
\noindent $({\cal A}_3)$ \parbox[t]{140mm}{$M\subset X$ is a neighborhood retract of $X$, 
$J_\alpha (N(M)) \subset M$ for $\alpha>0$, and  $M^* :=N(M)$ is an ${\cal L}$-retract (see \cite{Kryszewski-Ben-El} and \cite{Cw-Kr-JDE2009}), i.e.  there exist  a retraction $r:B(M^*,\eta)\to M^*$ with some $\eta>0$ and a constant $L>0$ such that
\be\label{17102011-1334}
\| r(\tau )- \tau\| \leq L d_{M^*} (\tau)  \ \mbox{ for all } \ \tau  \in B(M^*, \eta);
\ee
}\\[0.5em]
\noindent $({\cal A}_4)$ \parbox[t]{140mm}{  $F$ is continuous, bounded on bounded sets and satisfies the tangency condition
\be\label{17102011-1335}
F(N^{-1}(\tau))\in T_{M^{*}} (\tau),  \ \mbox{ for } \tau \in N(\overline U),
\ee
where $T_{M^*} (\tau)$ is the Bouligand tangent cone to $M^*$ at the point $\tau$, i.e.
$$
T_{M^*} (\tau): = \left\{ \theta \in X^* \mid \liminf_{\alpha \to 0^+} \frac{d_{M^*} (\tau + \alpha \theta) }{\alpha} =0 \right\}.
$$
}
\begin{Rem} \label{26032012-1932} {\em Since maximal monotone operators have closed graphs, it can be shown that in order to verify the continuity of the mapping ${\cal J}_{\ |\, X^* \times [\alpha_1, \alpha_2]}$ from condition $({\cal A}_2)$ it is sufficient to know that it maps bounded sets into relatively compact ones.}
\end{Rem}
Our goal is to transform the problem (\ref{17102011-1420}) into a fixed point one in $M$ and for which  fixed point index theory can be used. To this end define $\Phi_{\alpha}=\Phi_{\alpha}^{A,F}:\overline{U} \to M$ by
$$
\Phi_\alpha (u):= J_{\alpha} \left(r \left(N(u) + \alpha F(u)\right)\right), \ \ u\in \overline{U},  
$$
whenever $0<\alpha < \eta/\sup\{ \|F(u)\| \mid u\in \overline{U} \}$. Obviously, it is well defined, since for such $\alpha$ one has $(N+\alpha F) (\overline U) \subset B(M^*,\eta)$.
Moreover, observe that due to the assumptions, the mapping $r \circ (N+\alpha F)$ is bounded on bounded sets and, 
by $({\cal A}_2)$, $\Phi_\alpha$ is compact.\\
\indent Exploiting the tangency condition (\ref{17102011-1335}) and the inequality (\ref{17102011-1334}) together with compactness, we obtain the following localization of fixed points results.

\begin{Prop}\label{17102011-1441} If $K\subset \overline{U}$ is a closed set such that
$$
\{ u\in \overline U\cap D(A) \mid 0\in -Au + F(u)\} \cap K = \emptyset,
$$
then, for sufficiently small $\alpha>0$,
$\{ u\in \overline U \mid \Phi_\alpha (u)=u\} \cap K = \emptyset.$
\end{Prop}
\begin{Rem}\label{25102011-1802}{\em
Actually the tangency condition (\ref{17102011-1335}) and the continuity of $F\circ N^{-1}$ imply
$$
F\left( N^{-1}(\tau) \right) \in C_{M^*} (\tau) := \left\{ \theta \in X^* \mid \lim_{\alpha \to 0^+, \, \varrho \to \tau, \, \varrho\in M} \frac{d_M (\varrho+ \alpha \theta)}{\alpha} =0 \right\} \, \mbox{  for all } \tau \in N(U).
$$
Indeed
$$
F \left( N^{-1}(\tau) \right)=  
\lim_{\varrho \to \tau} F \left( N^{-1}(\varrho) \right) \in \Liminf_{\varrho\to \tau, \, \varrho \in M^*} T_{M^*} (\varrho) \subset C_{M^*} (\tau).
$$
The proof of the latter inclusion can be found in \cite{Aubin-Frankowska}.
}\end{Rem}
\begin{Lem}\label{17102011-1949}
{\em (i)} The graph $\Gr (A):= \{(u, \tau)\in X\times X^*  \mid u\in D(A)\}$ is closed;\\
{\em (ii)} If a sequence  of pairs $(u_n ,\tau_n) \in \Gr(A)$, $n\geq 1$, is bounded, 
then the sequence $(u_n)$ has a convergent subsequence.
\end{Lem}
{\bf Proof:} (i) Take any sequence of points $(u_n, \tau_n)\in \Gr(A)$, $n\geq 1$, such that 
$(u_n, \tau_n)\to (u_0, \tau_0)$ in $X\times X^*$, as $n\to +\infty$, for some $(u_0, \tau_0) \in X\times X^*$.
Clearly, $\tau_n \in Au_n$, and this gives $N(u_n) + \tau_n \in (N+A)(u_n)$, which, by $({\cal A}_1)$, gives $u_n = J_{1} (N(u_n)+\tau_n)$, $n\geq 1$. Hence, using the continuity of $N$ and $J_{1}$ yields 
$u_n = J_{1}(N(u_n)+\tau_n)) \to J_{1} (N(u_0) + \tau_0)$ as $n\to+\infty$, which implies 
$u_0 = J_{1} (N(u_0) + \tau_0)$, i.e. $\tau_0\in Au_0$. This shows that $\Gr(A)$ is closed.\\
\indent (ii) Note that, for each $n\geq 1$, $u_n = J_{1} (N(u_n) + \tau_n) \in J_{1} (N(B(0,R))+B(0,R))$, where $R>0$ is a constant such that $\|u_n\|_X \leq R$ and $\|\tau_n\|_{X^*} \leq R$ for $n\geq 1$. The boundedness of $N$ and $({\cal A}_2)$ imply that the set $(u_n)$ is a sequence of elements of the relatively compact set 
$J_{1} (N(B(0,R))+B(0,R))$. \hfill $\square$

\noindent {\bf Proof of Proposition \ref{17102011-1441}:} 
Suppose to the contrary that there exists a sequence $(\alpha_n)$ such that $\alpha_n\to 0^+$ such that for each 
$n\geq 1$ there is $u_n\in K$ with $\Phi_{\alpha_n} (u_n) = u_n$, that is
$$
N(u_n)+ \alpha_n \tau_n =  r(N(u_n) + \alpha_n F(u_n)) \mbox{ for some } \tau_n \in Au_n.
$$
In view of (\ref{17102011-1334}), one has
\begin{eqnarray}\label{26032012-1821}
\alpha_n \|\tau_n - F(u_n)\| &   =  & \| r(N(u_n) + \alpha_n F(u_n)) - ( N (u_n) + \alpha_n F (u_n) ) \| \\
                             & \leq & L d_{M^*} (N(u_n) + \alpha_n F(u_n)) \ \ \mbox{ for all } n\geq 1. \nonumber 
\end{eqnarray}
This implies 
$$
\| \tau_n \| \leq \|F(u_n)\| + L \alpha_{n}^{-1} d_{M^*} (N(u_n) + \alpha_n F(u_n)) \leq (1+L) \|F(u_n)\|, \ n\geq 1, 
$$
which means that $(\tau_n)$ is bounded. Therefore, by use of Lemma \ref{17102011-1949} (ii),
we may assume without loss of generality that $u_n \to u_0$ for some $u_0\in M$. 
Now using (\ref{26032012-1821})  and putting $p_n:=N(u_n)$, $n\geq 0$, we see that
$$
\|\tau_n - F(u_n)\| \leq L \cdot \frac{d_{M^*} (p_n + \alpha_n F(N^{-1}(p_0))}{\alpha_n} + L\|F(u_n)-F(u_0)\|, \ \ \mbox{ for } n\geq 1.
$$
By the tangency condition $({\cal A}_4)$ and Remark \ref{25102011-1802} together with the continuity of $F$, we get that $\tau_n \to F(u_n)$
as $n\to+\infty$. Hence, we have obtained that $(u_n, \tau_n)\to (u_0, F(u_0))$ and, by Lemma \ref{17102011-1949} (i),
$(u_0,F(u_0)) \in \Gr(A)$, i.e. $F(u_0)\in A u_0$, a contradiction completing the proof. \hfill $\square$

Now we put
\be\label{17102011-1439}   
\Deg_M ( A, F, U ) := \lim_{\alpha \to 0^+} \ind_M ( \Phi_\alpha , U ),    
\ee
where $\ind_M$ stands for the fixed point index 
for compact mappings of absolute neighborhood  retracts due to Granas -- see \cite{Granas} or \cite{Dugundji-Granas} for details.
We call this number as {\em the topological degree of coincidence  ({\em or just } topological degree) of $A$ and $F$ with respect to $M$}.

\begin{Th}\label{16122011-1102}
The coincidence degree defined by {\em (\ref{17102011-1439})} is well defined and has the following properties:\\
{\em (i)} \parbox[t]{140mm}{{\em (existence)} if $\Deg_M (A,F, U)\neq 0$, then there exists $u\in U\cap D(A)$ such that 
$0\in -Au+F(u)$;}\\[0.5em]
{\em (ii)} \parbox[t]{140mm}{{\em (additivity)} if $U_1, U_2$ are open disjoint subsets of a bounded open $U\subset M$ and $0\not\in (-A+F)(\o U \setminus (U_1\cup U_2))$, then
$$
\Deg_M (A,F, U) = \Deg_M (A,F, U_1) +\Deg_M (A,F, U_2);
$$}\\[0.5em]
{\em (iii)} \parbox[t]{140mm}{ {\em (homotopy invariance)}  if $H:\o U \times [0,1]\to X^*$ is a continuous and bounded mapping such that
$$
H( N^{-1}(\tau), t) \in T_{M^*} (\tau) \mbox{ for all } \tau \in N(\overline U), \, t>0,
$$
and $0\not\in -Au+H(u,t)$ for all $u\in \part U\cap D(A)$ and $t\in [0,1]$, then 
$$
\Deg_M (A,H(0,\cdot), U) = \Deg_M (A, H(1,\cdot), U);
$$}\\[0.5mm]
{\em (iv)} \parbox[t]{140mm}{{\em (normalization)} if $M$ is bounded and the mapping $\widetilde {\cal J}: X^* \times [0,+\infty) \ni (\tau, \alpha) \mapsto J^{\alpha} \tau \in X$ with $J^0 = N^{-1}$ is continuous, then $\Deg_M (A,F, M) = \chi (M)$.}
\end{Th}
{\bf Proof:} Note that for sufficiently small $\alpha>0$ it follows from Propostion \ref{17102011-1441}
that $\Phi_\alpha$ has no fixed point in $\part U$, i.e. the fixed point index $\ind_M (\Phi_\alpha, U)$ is well defined.
If $\alpha_1, \alpha_2>0$ are small enough, then, by $({\cal A}_2)$, $\Phi_{\alpha_1}$ is homotopic with $\Phi_{\alpha_2}$, which gives $\ind_M (\Phi_{\alpha_1}, U)= \ind_M (\Phi_{\alpha_2}, U)$, which means that 
the limit in (\ref{17102011-1439}) exists.\\
\indent (i) Suppose to the contrary that there is no $u\in U \cap D(A)$ such that 
$0\in - Au + F(u)$. Then, in view of Proposition \ref{17102011-1441}, for sufficiently small $\alpha>0$
the mappings $\Phi_\alpha$ have no fixed points in $\o U$, i.e. $\Deg_M (A,F,U) = \ind_M (\Phi_\alpha, U)=0$, a contradiction.\\
\indent  (ii) Due to Proposition \ref{17102011-1441}, for sufficiently small $\alpha>0$, $\Phi_\alpha$ has no fixed points
in $\o U \setminus (U_1\cup U_2)$. Therefore, by the definition of the degree,
$$
\Deg_M (A, F, U)=\ind_M (\Phi_\alpha, U) \ \  \mbox{ and } \ \ \Deg_M (A, F, U_k)=\ind_M (\Phi_\alpha, U_k) \mbox{ for } k=1,2.
$$
By the additivity property of the fixed point index
$$
\ind_M (\Phi_\alpha, U) = \ind_M (\Phi_\alpha, U_1) + \ind_M (\Phi_\alpha, U_2),
$$
which together with the earlier equalities gives the desired additivity of the degree.\\
\indent (iii) For sufficiently small $\alpha>0$ one can define $\Phi_\alpha:\o U\times [0,1] \to M$ by
$$
\Phi_\alpha (u,t):= J_{\alpha}\left( r (N(u) + \alpha H(u,t)) \right), \ u\in \o U, \ t\in [0,1]. 
$$
Proceeding along the lines of the proof of Proposition \ref{17102011-1441} we can prove that for sufficiently small 
$\alpha>0$
$$
\Phi_\alpha (u,t)\neq u \mbox{ for all } u\in \part U, \  t\in [0,1].
$$
Hence, by the homotopy invariance of the fixed point index and the formula defining the degree,
$$
\Deg_M (A,H(\cdot, 0),U) = \ind_M (\Phi_\alpha (\cdot,0), U) = \ind_M (\Phi_\alpha (\cdot,1), U)
=  \Deg_M (A,H(\cdot, 1),U).
$$
\indent (iv) Take small $\alpha>0$ such that $\Phi_\alpha$ is well defined. Then
$$
\Deg_M  (A, F, M) = \ind_M (\Phi_\alpha, M).
$$
Note that the normalization property for the fixed point index states that
the homomorphism $H_* (\Phi_\alpha):H_* (M)\to H_*(M)$ induced on (singular) homology spaces is a Leray endomorphism 
and
\be\label{19122011-1353}
\ind_M (\Phi_\alpha, M) = \Lambda (\Phi_\alpha)
\ee
where $\Lambda(\Phi_\alpha)$ is the generalized Leschetz number of the compact map $\Phi_\alpha$
-- see \cite[Definition V.(2.1), (3.1) and Theorem (5.1)]{Dugundji-Granas} or \cite{Granas}.
Further, consider $\Psi:M\times [0,1]\to M$ given by
$$
\Psi (u,t) := \widetilde {\cal J} (r(N(u)+t\alpha F(u)),t\alpha), \ \ u\in M, \ \ t\in [0,1].
$$
By the assumption, $\Psi$ is a continuous homotopy joining $\Psi(\cdot, 1) = \Phi_\alpha$
with the identity map $\id_M: M\to M$. Hence, for the maps induced on homology spaced one has 
$H_* (\Phi_\alpha) = H_* (\id_M)=\id_{H_* (M)}$ and, since $H_* (\Phi_\alpha)$ is an endomorphism Leray, we infer that
$\Lambda(\Phi_\alpha)=\sum_{n=0}^{\infty} (-1)^n \dim H_n (M) = \chi(M)$, which together with (\ref{19122011-1353}) ends the proof. \hfill $\square$

We end this section with a general result, which allows us to compute the degree is specific situations 
(comp. \cite[Prop. 4.2]{Cw-Kr-JDE2009}).
\begin{Th}\label{02122011-1318}
Let $M$ and $M^*$ be closed convex cones and that the mappings $A$ and $N$ are homogeneous with the same degree 
{\em (\footnote{i.e. there exists $\gamma>0$ such that $A(a u) = a^{\gamma} A(u)$, $u\in D(A)$, $a>0$, 
and $N(au)=a^\gamma N(u)$ for all $u\in X$, $a>0$.})}.
Suppose that there exists $\lma_1 \geq 0$ satisfying the following conditions\\[1em]
$({\cal M}_1)$ \parbox[t]{138mm}{$(A-\lma N)^{-1}(\{0\}) \cap  M = \{ 0\} \ \mbox{ for } \ \lma \neq \lma_1;$}\\[1em]
$({\cal M}_2)$ \parbox[t]{138mm}{there exists $\tau_0 \in M^*$ such that $(A-\lma N)^{-1} (\{ \tau_0 \}) \cap M =
\emptyset$ for $\lma>\lma_1$.}\\[1em]
Then
$$
\Deg_M (A, \lma N, B_M (0,\delta)) = \left\{ \begin{array}{ll} 1, & \lma < \lma_1,\\
0, & \lma>\lma_1,
\end{array} \right.
$$
for any $\delta>0$.
\end{Th}
{\bf Proof:}  Note that in view of $({\cal M}_1)$ the topological degree $\Deg_M (A, \lma N, B_M (0, \delta))$
is well defined.\\
\indent Now fix $\lma <\lma_1$. By the very construction, for sufficiently small $\alpha>0$,
\be\label{16122011-1528}
\Deg_M (A, \lma N, B_M (0,\delta)) = \ind_M (\Phi_\alpha, B_M (0,\delta))
\ee
where $\Phi_\alpha:\overline {B_M (0,\delta)}\to M$ is given by
$$
\Phi_\alpha (u):= J_\alpha (r(N(u) + \alpha \lma N(u))), \ u \in \overline {B_M (0,\delta)}.
$$
Define $\Theta: \overline {B_M (0,\delta)}\times [0,1]\to M$ by
$$
\Theta(u,t):= t \Phi_\alpha (u), \ u \in \overline {B_M (0,\delta)}, \, t\in [0,1].
$$
Suppose there are $u\neq 0$ and $t\in [0,1]$ such that $\Theta(u,t)=u$.
Then  $0 \in -A(u) + \mu N(u)$ with $\mu:=(t^\gamma-1)/\alpha + t^\gamma \lma$,  i.e. $u\in (A-\mu N )^{-1} (\{ 0\})\cap M$,
and, since $\mu=(t^\gamma-1)/\alpha + t^\gamma \lma < \lma_1$ we get a contradiction with $({\cal M}_1)$. Hence, we can use the homotopy
invariance of fixed point index to see that $\ind_M (\Phi_\alpha,B_M (0,\delta)) = 
\ind_M (0, B_M (0,\delta)) = 1$. This along with (\ref{16122011-1528}) implies the required equality.\\
\indent Let us pass to the case when $\lma>\lma_1$.  Define $H:M \times [0,1] \to X$ by
$H(u,t):= \lma N (u) + t\tau_0$, $u\in M$, $t\in [0,1]$. If $-A(u)+H(u,t)= 0$, then either $t=0$
and, due to $({\cal M}_1)$, $u=0$  or, by the homogeneity $-A(t^{-1/\gamma}u) + \lma N (t^{-1/\gamma}u) + \tau_0 = 0$,
where $\gamma>0$ is the common homogeneity degree for $A$ and $N$.
The latter equality contradicts $({\cal M}_2)$. Hence, the degrees $\Deg_M (A, H(\cdot,t), B_M (0,\delta) )$, $t\in [0,1]$,
are well defined and homotopy invariance can be used to obtain
$$
\Deg_M (A, \lma N, B_M (0,\delta)) = \Deg_M (A, \lma A +\tau_0, B_M (0,\delta).
$$
Finally the existence property of the degree together with $({\cal M}_2)$ implies $\Deg_M (A, \lma A +\tau_0, B_M (0,\delta))=0$, which completes the proof. \hfill $\square$

\section{Abstract setting for $p$-Laplacian} \label{29112011-1308}

Now we shall consider an abstract example falling into the setting of Section 2. 
It will be used in the sequel for the $p$-Laplace operator and the cone of positive functions in $L^{p}(\Omega)$.\\
\indent Let $X$ and $Y$ be reflexive normed spaces with a dense and compact linear embedding $i:Y\to X$.(\footnote{That is the mapping $i$ is linear and completely continuous with its range $i(Y)$ dense in $X$.})
Suppose that a closed convex cone $M\subset X$ and functionals ${\bold a}:Y\to \R$ and ${\bold n}: X\to \R$ satisfy the following conditions:\\[0.5em]
{\bf (a1)} \parbox[t]{62mm}{ ${\bold a}$ and ${\bold n}$ are coercive $C^1$ functionals;}(\footnote{By {\em coercivity} we mean that counterimages of intervals $(-\infty,m)$, with respect to a given functional, are bounded for all $m\in \R$.})\\[0.5em]
{\bf (a2)} \parbox[t]{140mm}{ there exists a continuous function $\kappa:[0,+\infty)\to [0, +\infty)$
such that $\kappa^{-1}(\{0\}) = \{0\}$, $\lim\limits_{s\to +\infty} \kappa (s) = +\infty$ and
\begin{eqnarray*}
\la D{\bold a} (u_1) - D{\bold a} (u_2), u_1-u_2 \ra_{Y} \geq \kappa (\|u_1-u_2\|_Y)\|u_1-u_2\|_Y \ \ \mbox{ for all } 
u_1, u_2\in Y,\\
\la D{\bold n} (u_1) - D{\bold n} (u_2), u_1-u_2 \ra_{X} \geq \kappa (\|u_1-u_2\|_X)\|u_1-u_2\|_X \ \ \mbox{ for all } 
u_1, u_2\in X;
\end{eqnarray*}
}\\[0.5em]
{\bf (a3)} \parbox[t]{140mm}{ for any $u\in M$ there exist $u^+, u^-\in M$ such that $u=u^+ - u^-$ and 
${\bold n}(u^+) \leq {\bold n}(u)$; if $u\in i(Y)$, then
$u^+, u^- \in i(Y)$ and ${\bold a} (i^{-1} u^+) \leq {\bold a} (i^{-1} u)$;
}\\[0.5em]
{\bf (a4)} \parbox[t]{140mm}{ $\bold n$ is bounded on bounded sets and monotone with respect to $M$, i.e. ${\bold n}(u+v) \geq  {\bold n} (u)$ for any $u,v\in M$.}\\[0.5em]

Let ${\cal A}:Y \to Y^*$ and $N:X\to X^*$ be defined by by $ {\cal A} : = D{\bold a}$ and $N:=D {\bold n}$. 
Note that that, due to ({\bf a2}), both ${\bold a}$ and ${\bold n}$ are strictly convex and ${\cal A}$ and $N$ are monotone operators. Define $A:D(A)\to X^*$ by 
\be\label{30112011-0001}
D(A):= i \left( {\cal A}^{-1} (i^*(X^*))\right) \mbox{ and } A u:= (i^*)^{-1}({\cal A} i^{-1}u), \mbox{ for } u\in D(A).
\ee
The above operation of restriction is a generalization of the analogical one that is usually considered in the case of a Gelfand triple $Y\subset X \subset Y^*$ where $X$ is a Hilbert space.\\

Below we show that assumptions $({\cal A}_1)$ and $({\cal A}_2)$ of Section 2 are satisfied.
\begin{Prop}\label{16122011-1117}
Under the above assumptions\\
{\em (i)} \parbox[t]{140mm}{  $N$ is a homeomorphism which is bounded on bounded sets;}\\[0.5em]
{\em (ii)} \parbox[t]{140mm}{ $N(M) = M^*:=\{ \tau\in X^* \mid \la \tau, u \ra \geq 0 \mbox{ for all } u\in M  \}$;}\\[0.5em]
{\em (iii)} \parbox[t]{140mm}{ $A$ is a densely defined maximal monotone operator;}\\[0.5em]
{\em (iv)} \parbox[t]{140mm}{ for any $\alpha>0$ and $\tau \in X^*$ there is a unique $u\in D(A)$ such that
$\tau = (N+\alpha A) (u)$;}\\[0.5em]
{\em (v)} \parbox[t]{140mm}{ if $J_\alpha: X^*\to X$, $\alpha >0$, is given by
$$
J_{\alpha} \tau := u \mbox{ where } u\in D(A) \mbox{ is such that } N(u)+\alpha A (u) = \tau,
$$
and ${\cal J}:X^* \times [0,+\infty) \to X$ by $${\cal J} (u, \alpha):= J_{\alpha} u,$$ 
then ${\cal J}$ is bounded on bounded sets and ${\cal J}_{\ |\, X^* \times [\alpha_1, \alpha_2]}$ with $0<\alpha_1 \leq \alpha_2$ is completely continuous;}\\[0.5em]
{\em (vi)} \parbox[t]{140mm}{$J_{\alpha}(M^*) \subset M$ for all $\alpha>0$.}
\end{Prop}
{\bf Proof:} To see (i), first note that $N$ is continuous, since ${\bold n}$ is $C^1$. Moreover, as a strictly convex
coercive functional on the reflexive Banach space $X$, for any $\tau\in X^{*}$, ${\bold n}-\tau$ admits a unique 
minimum point $u\in X$, i.e. $D{\bold n} (u)-\tau=0$, which gives  $N(u)=\tau$. Conversely, if $u\in X$ is such that
$N(u)=\tau$, then, by the strict convexity, $u$ is the unique minimum point. Hence, $N$ is bijective. To see that $N^{-1}$ is continuous, take any 
$(\tau_n)$ in $X$ with $\tau_n \to \tau$ in $X^*$ as $n\to +\infty$.
Observe that, by ({\bf a2}), we get
$$
\la \tau_n - \tau, N^{-1}(\tau_n)-N^{-1}(\tau)\ra_X  \geq \kappa (\|N^{-1}(\tau_n) - N^{-1} (\tau)\|_X)\|N^{-1}(\tau_n) - N^{-1} (\tau)\|_X,
$$
which yields the inequality
$$
\|\tau_n-\tau \|_{X^*}\geq \kappa\left(\|N^{-1}(\tau_n) - N^{-1} (\tau)\|_X \right).
$$
This in turn means that $N^{-1}(\tau_n) \to N^{-1} (\tau)$ in $X$ as $n\to +\infty$, that is $N^{-1}$ is continuous.
\indent To show that $N$ is bounded on bounded sets, we suppose to the contrary that there exists a bounded sequence $(u_n)$ in 
$X$ such that $\|N(u_n)\|_{X^*}\to +\infty$ as $n\to +\infty$. Since $X$ is reflexive, for each $n\geq 1$ one finds
an element $v_n\in X$ such that $\|v_n - u_n\|_{X} = 1$ and 
$$
\|N(u_n)\|_{X^*} = \la N(u_n), v_n-u_n\ra \leq {\bold n} ( v_n) - {\bold n} (u_n) \leq
\sup_{D_{X}(0, R+1)} {\bold n} - \inf_{D_X(0,R)} {\bold n} 
$$
where $R>0$ is such that $\|u_n\| \leq R$ for all $n\geq 1$. Thus, a contradiction proving the claim.\\
\indent To get (ii) take any $u\in M$ and $v\in M$. In view of  ({\bf a4})
$$
{\bold n}(u+h v)-{\bold n} (u) \geq 0 \  \mbox{ for any } h>0,
$$
which, after a division by $h$ and passage to the limit with $h\to 0^+$, yields $\la N(u), v\ra \geq 0$.
Hence $N(M)\subset M^*$. To prove the converse inclusion $M^*\subset N(M)$, we take any $\tau\in M^*$. 
As we mentioned ${\bold n}-\tau$ attains the minimum at some $u\in X$. On the other hand, by ({\bf a2}),
$$
{\bold n}(u^+) -\tau (u^+) \leq {\bold n}(u) - \tau(u^+) + \tau (u^-) = {\bold n}(u) - \tau(u).
$$ 
and, since the minimum point is unique, we infer that $u=u^+ \in M$.\\
\indent To show (iii), take any $u_1, u_2\in D(A)$. Clearly  $(A u_k) \circ i = {\cal A} (\tilde u_k)$ with $\tilde u_k = i^{-1} (u_k)$, for $k=1,2$.
Therefore, by ({\bf a2}),
$$
\la A u_1-Au_2, u_1 - u_2 \ra_{X} = [Au_1 - Au_2] i(\tilde u_1-\tilde u_2)
=  \la {\cal A} (\tilde u_1) - {\cal A} (\tilde u_2), \tilde u_1-\tilde u_2\ra_Y \geq 0.
$$
Hence $A$ is monotone and it is left to prove that $A$ is maximal monotone, i.e. that additionally one has $A (D(A))=X^*$.
To see it we choose any $\tau \in X^*$ and put $\Phi:={\bold a} - i^*(\tau)$. $\Phi$ 
is a convex coercive functional on the reflexive space $Y$. Hence it admits a miniumem, i.e. there is a point 
$\bar u\in Y$ such that $D\Phi(\bar u)=0$, i.e. ${\cal A}(\bar u) = i^*(\tau)$. This means that
$u:=i(\bar u)\in D(A)$ and that $A(u)=\tau$.\\
\indent To show (iv) take any $\tau \in X^*$ and $\alpha>0$. We proceed like in (iii), that is we consider a functional 
$\Phi:= {\bold n}\circ i + \alpha {\bold a} - i^*(\tau)$ on $Y$. It is clear that $\Phi$ -- as a strictly convex and coercive functional on a reflexive Banach space -- admits a minimum, i.e. there exists $\bar u\in Y$ such that 
$D\Phi(\bar u)=0$. This means that $i^*(N(i(\bar u)))  + \alpha {\cal A} (\bar u) = i^*(\tau)$. Subsequently, we deduce that ${\cal A}(\bar u) \in i^*(X^*)$, i.e. $u:= i(\bar u)\in  D(A)$ and $N(u) + \alpha A (u) = \tau$. Moreover, 
observe that for each $u\in D(A)$ such that $N(u) + \alpha A (u) = \tau$, $i^{-1}(u)$ is a critical point of $\Phi$. Since $\Phi$ is strictly convex it has to be the unique minimum point.\\
\indent (v) Suppose that a sequence $(\tau_n)$ is bounded in $X^*$ and $(\beta_n)$ is a sequence in $[\alpha_1, \alpha_2]$. Put $u_n:= J_{\beta_n} (\tau_n)$, $n\geq 1$.
Then $i^*N(u_n) + \beta_n {\cal A} (\bar u_n) = i^*(\tau_n)$, where $\bar u_n:=i^{-1}(u_n)$, $n\geq 1$.
Since  $N$ is bounded and $\beta_n > \alpha_1>0$ for all $n\geq 1$, we infer that  $({\cal A}(\bar u_n))$ is bounded. Observe that, in view of ({\bf a2}),
$$
\la {\cal A} (\bar u_n) - {\cal A}(0), \bar u_n \ra_{Y} \geq \kappa (\|\bar u_n\|_{Y})\|\bar u_n\|_Y,
$$
i.e. $\|{\cal A}(u_n) - {\cal A}(0)\|_{Y} \geq \kappa (\|\bar u_n\|_{Y})$. Hence, by the boundedness of $({\cal A}(\bar u_n))$ and the assumed property of $\kappa$, $(\bar u_n)$ is bounded. Therefore $(u_n) = (i(\bar u_n))$ is relatively compact, which together with Remark \ref{26032012-1932} proves the assertion. \\
\indent In order to prove (vi), take any $\tau \in M^*$. We need to show that $u:=J_{\alpha} (\tau) \in M$.
In the proof of (iv) we have showed that $i^{-1} u$ is the unique minimum of the functional 
$\Phi={\bold n}\circ i + \alpha {\bold a} - i^*(\tau)$ on $Y$. On the other hand, by use of ({\bf a3}) and the definition of $M^*$, one has
$$
\Phi(i^{-1} u^+ )\!=\! {\bold n} (u^+) \!+\! \alpha {\bold a}(i^{-1} u^+)\! -\! \tau (u^+)\!\leq\! 
{\bold n}(u) \!+\! \alpha {\bold a} (i^{-1} u))\!-\!\tau ( u^+) \!+\! \tau(u^-) = \Phi(i^{-1}u). 
$$ 
This means that $i^{-1} u=i^{-1} u^+$ and $u\in M$. \hfill $\square$

\section{Elliptic problems with $p$-Laplacian}

Now we shall apply the above abstract setting from the previous section to the $p$-Laplacian problem. To this end fix $p>2$, and put
$$
X_p:= L^p (\Omega), \ Y_p:=W_{0}^{1,p} (\Omega) \mbox{ and } M_p:=\{u\in X \mid u(x)\geq 0 \mbox{ for a.e. } x\in \Omega\}.
$$ 
Both, $X_p$ and $Y_p$ are reflexive and, by the Rellich-Kondrachov theorem, the natural embedding $i:Y_p\to X_p$ is compact and dense. It is easy to see that $M_p$ is a closed convex subset of $X_p$.
Next define functionals ${\bold a}:Y_p\to \R$ and ${\bold n}:X_p\to\R$ by
\begin{eqnarray*}
{\bold a} (u):= \frac{1}{p} \int_{\Omega} |\nabla u(x)|^p \d x, \ u\in Y_p,\\
{\bold n}(u):= \frac{1}{p} \int_{\Omega} |u(x)|^p \d x, \ u\in X_p.
\end{eqnarray*}
We prove that these objects satisfy the abstract assumptions of the general setting.
\begin{Prop} \label{16122011-1119}
The functionals ${\bold a}$ and ${\bold n}$ with the cone $M_p$ satisfy all the assumptions
{\em {\bf (a1) -- (a4)}} from Section {\em \ref{29112011-1308}} and
\begin{eqnarray}
\la D{\bold a} (u), v\ra_{Y} = \frac{1}{p} \int_{\Omega} |\nabla u(x)|^{p-2} \nabla u(x)\cdot \nabla v(x) \d x, \ u, v \in Y_p, \label{29112011-1313}\\
\la D{\bold a} (u) - D{\bold a} (v), u-v \ra_Y \geq 2^{2-p}\|u-v\|_{Y}^{p}, \, u,v\in Y_p, \label{30112011-2134}\\  
\label{29112011-1314} \la D{\bold n} (u), v\ra_{X} = \frac{1}{p} \int_{\Omega} |u(x)|^{p-2} u(x) v(x) \d x, \ u,v \in X_p,\\
\la D {\bold n} (u) - D {\bold n} (v), u-v \ra_X \geq 2^{2-p}\|u-v\|_{X}^{p}, \, u,v\in X_p.\label{30112011-2135}
\end{eqnarray}
Moreover, if $A_p:D(A_p) \to X_p$ is defined, in analogy to {\em (\ref{30112011-0001})}, by
$$
D(A_p):= i\left( (D {\bold a})^{-1} (i^*(X^*))\right) \mbox{ and } A_p u:= (i^*)^{-1}( D({\bold a}) i^{-1} u), \mbox{ for } u\in D(A),
$$
then
$$
A_p u= -\mathrm{div} (|\nabla u|^{p-2} \nabla u), \mbox{ for } u\in D(A_p),
$$
where the divergence is meant in the distributional sense and 
$$D(A_p)=\{u\in W_{0}^{1,p} (\Omega) \mid \mathrm{div}(|\nabla u|^{p-2}\nabla u) \mbox{ exists and belongs to }  L^{p} (\Omega)\}.$$
\end{Prop}
{\bf Proof:} In order to see {\bf (a1)}, note that the functionals ${\bold a}$ and ${\bold n}$ are clearly Gateaux differentiable with the formulas 
(\ref{29112011-1313}) and (\ref{29112011-1314}) satisfied. Since these Gateaux derivatives 
are continuous the functionals are Fr\^{e}chet differentiable. The coercivity is immediate 
as ${\bold a}(u) = (1/p)\|u\|_{Y_p}^{p}$, $u\in Y$, and ${\bold n}(u)=(1/p)\|u\|_{X_p}^{p}$, $u\in X$.\\
One can check the condition {\bf (a2)}, i.e. (\ref{30112011-2134}) and (\ref{30112011-2135}), by use of the following inequality
\be\label{10122011-1810}
(|x|^{p-2}x-|y|^{p-2}y)\cdot (x-y)\geq 2^{2-p} |x-y|^p \mbox{ for any  } x,y\in \R^M, \, M\geq 1.
\ee
Obviously, for $\kappa:[0,+\infty)\to [0, +\infty)$, given by $\kappa (s):=2^{2-p} s^p$, $s\geq 0$, 
one has $\kappa^{-1}(\{0\}) = \{0\}$, $\lim\limits_{s\to +\infty} \kappa (s) = +\infty$.\\
\indent As for {\bf (a3)}, take any $u\in X$. Then taking $u_+ := \max\{ u , 0\}$ and $u_-:=\max\{-u,0\}$
we have $u= u^+ - u^-$ and 
$$
{\bold n}(u_+ ) = \frac{1}{p} \int_{\Omega} |u_+ (x)|^p \d x \leq \frac{1}{p} \int_{\Omega} |u (x)|^p \d x  = {\bold n}(u).
$$
If $u\in Y_p = W_{0}^{1,p}(\Omega)$, then, due to Lemma 7.6 of \cite{Gilbarg-Trudinger},
$\nabla u_+(x) = 0$ if $u(x)\leq 0$ and $\nabla u_+ (x)=\nabla u(x)$ if $u(x)\geq 0$.
Therefore $u_+ \in Y_p$ and
$$
{\bold a}(u_+) = \frac{1}{p}\int_{\Omega}|\nabla u_+(x)|^p \d x \leq  \frac{1}{p}\int_{\Omega}|\nabla u (x)|^p \d x
= {\bold a}(u). 
$$
\indent Finally, {\bf (a4)} is immediate as, for $u,v\in M_p$, $|u|^{p}= u^p \leq (u+v)^{p}=|u+v|^p$. \hfill $\square$

In view of Section \ref{29112011-1308}, the operators $A_p$, $N_p:=D {\bold n}$ together with $M_p$ and $M_p^*$ satisfy the assumptions made in Section \ref{30112011-2213} and the topological degree can be applied for perturbations of $A_p$.
Before we proceed further let us pay attention to the perturbation term.
\begin{Prop}\label{30112011-2240}
Let $f:\Omega \times [0,+\infty)\to \R$ satisfy {\em (\ref{15112011-2220})} and $f(x,0) \geq 0$ for a.a. $x\in \Omega$.
Then the mapping $F : X_p \to X_{p}^{*}$ given by
$$
\la F(u), v \ra_{X_p}:= \int_{\Omega} f(u(x))v(x) \d x, \ \ u\in M_p, \, v\in X_p,
$$
is well defined, continuous, bounded on bounded sets and
\be\label{30112011-2346}
F(N^{-1}(\tau)) \in T_{M_{p}^*} (\tau) \ \ \mbox{ for  any } \ \tau\in M_p^{*}.
\ee
\end{Prop}
\begin{Lem}\label{30112011-2348}
Let $1<q<\infty$, $\Omega\subset \R^N$, $N\geq 1$, be open and
$$
M_q:= \{ u \in L^q(\Omega) \mid u(x) \geq 0 \mbox{ for a.e. } x\in \Omega\}.
$$
Then $\ \  T_{M_q} (u) = \{v \in L^q (\Omega) \mid v(x)\geq 0 \mbox{ for a.e. } x\in \Omega \mbox{ such that } u(x)=0 \}$.
\end{Lem}
{\bf Proof}: Put $T_u:=\{v \in L^q (\Omega) \mid v(x)\geq 0 \mbox{ for a.e. } x\in \Omega \mbox{ such that } u(x)=0 \}$.
To see that $T_u\subset T_{M_q} (u)$ take any $v\in T_u$ and define $v_n\in L^1(\Omega)$, $n\geq 1$, by
$$
v_n (x):=\left\{ \begin{array}{cl} 
v(x) & \mbox{ if } v(x) + n u(x) \geq 0,\\
0 & \mbox{ if } v(x) + n u (x) <0. \end{array}\right.
$$
Clearly, $v_n \in M_q - n u \subset T_{M_q} (u)$, for each $n\geq 1$.
Moreover it is clear that, for a.e. $x\in \Omega$ and any $n\geq 1$, $v_n(x) = v(x)\geq 0$ if $u(x)=0$ 
and $v_n (x) \to v(x)$ if $u(x)>0$. This implies that $v_n \to v$ in $L^q(\Omega)$, i.e. $v\in T_{M_q} (u)$.\\
\indent  In order to show the converse inclusion, observe that $T_u$ is closed and, for any $h>0$,
$h(M-u) \subset T_u$. This clearly implies that $T_{M_q} (u)\subset T_u$. \hfill $\square$

\noindent {\bf Proof of Proposition \ref{30112011-2240}}:
Using the Riesz representation isomorphism $\varrho$ between $L^p (\Omega)^*$ and $L^q(\Omega)$, $1/p +1/q=1$,
the mapping $F\circ N^{-1}$ can be treated as the mapping $L^q(\Omega)\ni u \mapsto f(\cdot, \theta_q (u)) \in L^q(\Omega)$	
where $\theta_q:\R\to\R$ is given by $\theta_q (s)=|s|^{q-2}s$, $s\in\R$. It is well defined as
$$
|(f(x,\theta_q(s))| \leq C(1+|\theta_q (s)|^p) = C(1+|s|)  \mbox{ for } s\geq 0 \mbox{ and a.e. } x\in\Omega.
$$
Observe that $f(x,\theta_q (0))=f(x,0)\geq 0$ for. a.e. $x\in\Omega$, which, by use of Lemma \ref{30112011-2348}, 
implies that $f(\cdot, \theta_q (u(\cdot)))\in T_{M_q}(u)$ for all $u\in M_q$. 
Since $\varrho (M_p^*)=M_q$, we infer that (\ref{30112011-2346}) holds. \hfill $\square$

Hence we have showed that the problem (\ref{29102011-1511}) indeed can be formulated as an abstract 
problem 
$$
\left\{\begin{array}{l}
A_p (u) = F (u),\\
u\in M_p \cap D(A_p). \end{array}
\right.
$$
In order to take advantage of the topological degree effectively we need some methods of computing it. 
\begin{Th}\label{02122011-1924}
If  $\ 2<p<\infty$ and $\rho\in L^{\infty} (\Omega)$ is such that either $\rho (x)> \lambda_{1,p}$ for a.e. $x\in \Omega$,
or $\rho(x)< \lambda_{1,p}$ for a.e. $x\in \Omega$, then
$$
\Deg_{M_p} (A_p, \rho N_{p}, B_{M_p} (0,R)) = \left\{ 
\begin{array}{cl}
1, & \mbox{ if } \ \ \rho (x) < \lambda_{1,p} \ \  \mbox{ for a.e. } \ \ x\in \Omega,\\
0, & \mbox{ if } \ \ \rho(x) > \lambda_{1,p} \ \ \mbox{ for a.e. } \ \ x\in \Omega.
\end{array}
\right.
$$
\end{Th}
\begin{Rem}\label{04122011-1324} 
{\em Before passing to the proof of Theorem \ref{02122011-1924}, we need to make a comment on the 
eigenvalue problem relating to  the $p$-Laplace operator. 
Solving the nonlinear eigenvalue problem 
$$
\left\{\begin{array}{l}
A_p (u) = \lma N_p(u)\\
u\in M_p \cap D(A_p)
\end{array} \right.
$$
reduces to find nonnegative weak solutions $u\in W^{1,p}(\Omega)$ of 
\be\label{02122011-2016}
\left\{ \begin{array}{cl}
-\mathrm{div} (|\nabla u(x)|^{p-2} \nabla u(x)) = \lma |u(x)|^{p-2} u(x), & x\in \Omega,\\
u(x) = 0, & x\in \partial \Omega. \\
\end{array}\right. 
\ee
It appears that some properties of the eigenvalue problem for the Laplace operator are also valid 
for the $p$-Laplace one. For details we refer to \cite{Lindqvist}, \cite{Lindqvist-PAMS-1992} and 
\cite{Lindqvist-1}. In particular, it is known that (\ref{02122011-2016}) does not admit any nonzero solutions if $\lambda\leq 0$, i.e. the $p$-Laplace has no nonpositive eigenvalues. Moreover, there exists the smallest eigenvalue $\lma_{1,p}$ given by the Rayleigh formula
$$
\lambda_{1,p} = \inf\limits_{u\in W_{0}^{1,p} (\Omega), u\neq 0} \frac{\int_{\Omega} |\nabla u(x)|^p \d x}{\int_{\Omega} |u(x)|^p \d x}.
$$
The eigenfunctions corresponding to $\lma_{1,p}$ are either strictly positive or negative in $\Omega$ and belong to 
$L^\infty (\Omega)$. Moreover, $\lambda_{1,p}$ is an isolated eigenvalue and if there are two eigenfunctions $u,v$ for $\lambda_{1,p}$, then there exists $\alpha\in\R$ such that $u=\alpha v$. It is also known that if any eigenfunction does not change its sign in 
$\Omega$, then the corresponding eigenvalue must be equal to $\lambda_{1,p}$. \hfill $\square$}
\end{Rem}
In the proof we shall use  a few lemmata given below.
\begin{Lem}\label{10122011-1907}
There are $ C, s > 0$ such that $\|u\|_{L^p} \leq C |\tilde \Omega|^{s} \|\nabla u\|_{L^p}$ for
all $u\in W_{0}^{1,p}(\Omega)$ and measurable $\tilde \Omega\subset \Omega$ with the property
$u(x)=0$ if $x\not\in \tilde \Omega$.
\end{Lem}
{\bf Proof:} By the Sobolev embedding theorem there exists $q>p$ such that 
$$
\| u \|_{L^q} \leq C\|\nabla u\|_{p} \mbox{ for all } u\in W_{0}^{1,p}(\Omega).
$$ 
On the other hand, by the H\"{o}lder inequality,
$$
\|u\|_{L^p} \leq \| u \|_{L^q}   |\tilde \Omega|^{1/p-1/q}. 
$$
Combining the two above inequalities we get the desired one with $s:=1/p-1/q$. \hfill $\square$

\begin{Lem}\label{09122011-1525}
Let $v$ be a nonnegative weak solution of {\em (\ref{02122011-2016})} with $\lma=\lma_{1,p}$ and $\rho\in L^{\infty}(\Omega)$. If $u\in W_{0}^{1,p}(\Omega)$ is a weak solution to 
$$
-\mathrm{div} (|\nabla u|^{p-2} \nabla u) = \rho |u|^{p-2} u  + |v|^{p-2}v \, \mbox{ on } \Omega,
$$
then $u\in L^{\infty}(\Omega)$.
\end{Lem}
{\bf Proof}: Here we adapt the arguments from \cite{Lindqvist-PAMS-1992}. 	Note that without loss of generality
we can consider the equation
$$
-\mathrm{div} (|\nabla u|^{p-2} \nabla u) = \rho |u|^{p-2} u  + \lambda_{1,p} |v|^{p-2} v  \, \mbox{ on } \Omega.
$$
Take any $k>0$ and put $\eta:=\max\{u-v-k,0\}$. Since $\eta\in W_{0}^{1,p}(\Omega)$, we get
$$
\int_{\Omega_k} (|\nabla u|^{p-2}\nabla u - |\nabla v|^{p-2} \nabla v) \cdot \nabla (u-v)\d x \leq \|\rho\|_{L^{\infty}} \int_{\Omega_k} u^{p-1}(u-v-k) \d x
$$
with $\Omega_k := \{x\in \Omega \mid u(x)-v(x)-k > 0 \}$. This, by use of (\ref{10122011-1810}) and the convexity
of the function $s\mapsto |s |^{p-1}$, gives
\begin{eqnarray*}
\int_{\Omega_k} |\nabla (u-v)|^p \d x  & \leq & C_1 \int_{\Omega_k} u^{p-1}  (u-v-k) \d x \\
& \leq &  C_1 2^{p-2} \left( \int_{\Omega_k} (u-v-k)^p \d x + \int_{\Omega_k} (v+k)^{p-1}(u-v-k) \d x\right)
\end{eqnarray*}
for some constant $C_1>0$ (here all the constant are to be independent of $k$). 
By applying Lemma \ref{10122011-1907}, one gets 
$$
\int_{\Omega_k} (u-v-k)^p \d x \leq C|\Omega_k|^{s} \int_{\Omega_k} |\nabla (u-v)|^{p} \d x,
$$
which together with the previous inequality yields
$$
(1- C_2|\Omega_k|^{s}) \int_{\Omega_k} (u-v-k)^p \d x \leq C_2 |\Omega_k|^{s} \int_{\Omega_k} (v+k)^{p-1}(u-v-k) \d x
$$
for some $C_2>0$. Since $|\Omega_k|\to 0$ as $k\to +\infty$, there is $k_0$ such that for all $k\geq k_0$
$1-C_2 |\Omega_k|^{s} >1/2$. Further, for $k\geq k_0$,
$$
 \int_{\Omega_k} (u-v-k)^p \d x   \leq  2 C_2 |\Omega_k|^{s} (\|v\|_{L^{\infty}} + k)^{p-1} \int_{\Omega_k} (u-v-k) \d x.
$$
Next we observe that the H\"{o}lder inequality yields
\be\label{10122011-2019}
\int_{\Omega_k} (u-v-k) \d x \leq C_4 k |\Omega_k|^{1+s(p-1)^{-1}}  \mbox{ for all } k\geq k_0
\ee
and some constant $C_4>0$. Now define  $j:(0,+\infty)\to [0,+\infty)$ by
$$
j(k):= \int_{\Omega_k} (u-v-k) \d x, \, \, k>0.
$$
Note that by the Tonelli-Fubini theorem applied to the set $\{(x,t)\in \Omega\times [0,+\infty)\mid u(x)-v(x)>t>k \}$
one has
$$
j(k) = \int_{k}^{+\infty} |\Omega_t| \d t, \,\,  k> 0.
$$
Obviously, $j$ is nonincreasing and absolutely continuous with $j'(k)=-|\Omega_k|$ for a.e. $k\geq 0$.
We claim that $j(k)=0$ for some $k>0$. If it were not so, then (\ref{10122011-2019}) could  be rewritten as 
$$
j(k)^{\theta} \leq -C_{4}^{\theta} k^{\theta} j'(k) \mbox{ for all } k\geq k_0
$$
with $\theta:=(1+s(p-1)^{-1})^{-1}$, and consequently
$$
k^{-\theta} \leq - C_{4}^{\theta} j(k)^{-\theta} j'(k) \mbox { for all } k\geq k_0.
$$
This after integration would give
$$
k^{1-\theta} + C_{4}^{\theta} j(k)^{1-\theta} \leq k_0^{1-\theta} + C_{4}^{\theta} j(k_0)^{1-\theta} \mbox{ for all } k\geq k_0,
$$
which yields a contradiction proving the claim that $j(k)=0$ for some $k>0$.
Then, for some $k>0$, $|\Omega_k|=0$ and $u\leq v+k$ a.e. on $\Omega$. This shows that 
$u\in L^{\infty}(\Omega)$, as 
$v\in L^{\infty}(\Omega)$ (see Remark \ref{04122011-1324}). \hfill $\square$

\begin{Lem}{\em (see \cite[Th. 1]{Fleckinger-et-al})}\label{02122011-2037}
If $h\in  L^\infty(\Omega)$ is nonnegative and nonzero, then the equation
$$
-\mathrm{div} (|\nabla u|^{p-2} \nabla u) = \lambda_{1,p} |u|^{p-2}u + h, \, \mbox{ on } \Omega,
$$ 
has no nonzero weak solution in $W_{0}^{1,p}(\Omega)$.
\end{Lem}
\begin{Lem}\label{11122011-0003}
If $\rho\in L^\infty (\Omega)$ and either $\rho (x)> \lambda_{1,p}$ for a.e. $x\in\Omega$ or $\rho(x)<\lambda_{1,p}$ for a.e. $x\in \Omega$, 
then the problem 
\be\label{10122011-2336}
-\mathrm{div}(|\nabla u|^{p-2}\nabla u) = \rho |u|^{p-2}u \mbox{ on } \Omega
\ee
does not admit a nonzero solution $u\in W_{0}^{1,p}(\Omega)$ such that $u\geq 0$.
\end{Lem}
{\bf Proof}:  If $\rho<\lambda_{1,p}$ a.e. on $\Omega$ and $u\in W_{0}^{1,p}(\Omega)$ is a nonzero weak solution of (\ref{10122011-2336}), then
$$
\int_{\Omega} |\nabla u|^{p} \d x = \int_{\Omega} \rho |u|^{p-2}u \d x < \lambda_{1,p} \int_{\Omega} |u|^{p-2}u \d x,
$$
which gives $\lma_{1,p}> \int_{\Omega} |\nabla u|^p \d x / \int_{\Omega} |u|^p \d x$, a contradiction with the Rayleigh formula.\\
\indent  In the case  $\rho>\lambda_{1,p}$ a.e. on $\Omega$, we observe that if $u$ is a weak solution of (\ref{10122011-2336}), then $u$ is a weak solution of
$$
-\mathrm{div}(|\nabla u|^{p-2}\nabla u) =\lma_{1,p} |u|^{p-2}u + h  \ \mbox{ on } \ \Omega
$$  
with $h:=(\rho - \lambda_{1,p})|u|^{p-2}u$. Clearly, $h\geq 0$ and $h\in L^{\infty}(\Omega)$, since 
$u\in L^{\infty}(\Omega)$ due to Lemma \ref{09122011-1525}. Hence, Lemma \ref{02122011-2037} leads to a contradiction ending the proof. \hfill $\square$

\noindent {\bf Proof of Theorem \ref{02122011-1924}:} 
Assume that $\rho>\lambda_{1,p}$ a.e. on $\Omega$ and fix $\tilde \lambda>\lambda_{1,p}$.
Define $H:X_p \times [0,1]\to X_p$ by $H (u,t) := (t\tilde \lambda + (1-t) \rho) N_p(u)$, $u\in X_p$, $t\in [0,1]$.
In view of Lemma \ref{11122011-0003}, $-A_p (u)+H(u,t)\neq 0$ for all $u\in D(A_p) \setminus \{ 0 \}$ and $t\in [0,1]$.
Therefore, we can use the homotopy invariance 
-- Theorem \ref{16122011-1102} (iii) to get
\be\label{11122011-0122}
\Deg_{M_p} (A_p, \rho N_{p}, B_{M_p} (0,R) ) = \Deg_{M_p} ( A_p, \tilde \lma N_p, B_{M_p}(0,R) ). 
\ee
In a similar manner one can prove the same formula in the case $\rho<\lambda_{1,p}$ a.e. on $\Omega$ with $\tilde \lambda
<\lambda_{1,p}$.\\
\indent Now we shall prove that conditions $({\cal M}_1)$ and $({\cal M}_2)$ of Theorem \ref{02122011-1318} are satisfied.
Observe that, in view of Lemma \ref{11122011-0003}, for any $\lma \neq \lma_{1,p}$, the eigenvalue problem (\ref{02122011-2016}) has no nontrivial and nonnegative weak solutions, i.e. $({\cal M}_1)$ holds. 
To show $({\cal M}_2)$ let $\tau_0\in L^p(\Omega)$ be the functional determined by $|u_0|^{p-2}u_0$ with $u_0$
being a fixed positive solution of the eigenvalue problem (\ref{02122011-2016}) with $\lambda=\lambda_{1,p}$.
Suppose that there exists $u\in (A_p-\lambda N_p)^{-1}( \{ \tau_0 \} ) \cap M_p$ for some $\lambda > \lambda_{1,p}$.
This means that $u\in W_{0}^{1,p}(\Omega)$ is a nonnegative weak solution of
$$
-\mathrm{div} (|\nabla u|^{p-2} \nabla u)  = \lambda_{1,p} |u|^{p-2} u + h \ \mbox{ on } \ \Omega
$$
with $h:=(\lambda-\lambda_{1,p}) |u|^{p-2} u + |u_0|^{p-2}u_0$.  It follows from Lemma \ref{09122011-1525} that $h\in L^{\infty} (\Omega)$. Since $h\geq 0$, Lemma \ref{02122011-2037} implies that such a solution does not exist, a contradiction proving $({\cal M}_2)$. Hence, by Theorem \ref{02122011-1318} and (\ref{11122011-0122}), the desired formula follows. \hfill $\square$

\noindent The obtained formula  results in the following general one. 
\begin{Th}\label{16122011-1341}
Let $f$ and $F$ be as in Proposition \ref{30112011-2240} and suppose that {\em (\ref{15112011-2220})} hold.\\
\indent {\em (i)} \parbox[t]{135mm}{If $\rho_0$ is as in {\em (\ref{15112011-2222})} and either 
$\rho_0 (x) < \lambda_{1,p}$, for a.e. $x\in\Omega$, or $\lambda_{1,p}<\rho_0 (x)$, for a.e. $x\in\Omega$, 
then there exists $\delta>0$ such that
$A_p (u) \neq F (u) \mbox{ for all } u \in D(A_p) \cap \left( B_{M_p} (0,\delta) \setminus \{ 0 \}\right)$ and
$$
\Deg_M (A, F, B_M (0,\delta)) = \left\{ 
\begin{array}{cl}
1, &  \ \ \mbox{ if } \ \ \rho_0(x) < \lambda_{1,p} \ \ \mbox{ for a.e. } x\in \Omega,\\
0, & \ \ \mbox{ if } \ \ \rho_0(x) > \lambda_{1,p} \ \ \mbox{ for a.e. } x\in \Omega.
\end{array}
\right.
$$}\\[0.5em]
\indent {\em (ii)} \parbox[t]{135mm}{ If $\rho_\infty$ is as in {\em (\ref{15112011-2222})} 
either  $\rho_\infty (x) < \lambda_{1,p}$, for a.e. $x\in\Omega$, or $\lambda_{1,p}<\rho_\infty (x)$, for a.e. $x\in\Omega$,
then there exists $R>0$ such that
$A_p (u) \neq F(u) \mbox{ for all } u \in D(A_p) \cap \left( M_p\setminus B_{M_p} (0,R)\right)$  and
$$
\Deg_{M_p} (A_p, F, B_{M_p} (0,R)) = \left\{ 
\begin{array}{cl}
1, & \mbox{ if }  \ \ \rho_\infty (x)< \lambda_{1,p} \ \ \mbox{ for a.e. } x\in \Omega,\\
0, & \mbox{ if }  \  \ \rho_\infty (x)> \lambda_{1,p} \ \ \mbox{ for a.e. } x\in \Omega.
\end{array}
\right.
$$}
\end{Th}
{\bf Proof}: (i) Define $H:M_p\times [0,1]\to X_p$ by 
$H(u,t):= t F(u) + (1-t) \rho_0 N_{p} (u),$ $(u,t)\in M_p \times [0,1]$.
By Proposition \ref{30112011-2240}, $H$ is continuous and $F\circ N_{p}^{-1}$ is tangent to $M^*$. 
Moreover we claim that 
\be\label{16122011-1303}
\mbox{ there is $\delta >0$ 
such that $-A_p (u)+H(u,t) \neq 0$ for all $u\in M_p \cap D(A_p),  \ t\in [0,1]$.}
\ee
Suppose to the contrary that there exists $(u_n)$ in $(M_p\cap D(A_p))\setminus \{ 0 \}$ and $(t_n)$ in $[0,1]$ such that
$u_n \to 0$ in $X_p$ and $-A_p (u_n) + H(u_n, t_n)=0$, $n\geq 1$. Then clearly, if we put $w_n:= \|u_n\|_{X_p}^{-1}u_n$
and $s_n:=\|u_n\|_{X_p}$, then
$A_p(w_n) = s_{n}^{1-p} H(s_n w_n, t_n)$, which gives
\be\label{16122011-1248}
w_n  =J_{1} \left( N_p (w_n) + s_{n}^{1-p} H(s_n w_n, t_n) \right), \ n\geq 1.
\ee
The growth condition (\ref{15112011-2220}) and the existence of the first limit in (\ref{15112011-2222}) imply that there exists $C_1>0$ such that $\|N_p (w_n) + s_{n}^{1-p} H(s_n w_n, t_n)\|_{X_{p}^{*}}\leq C_1$ for all $n\geq 1$.
Therefore we infer that $(w_n)$ has a subsequence convergent in $X_p$,
since, according to Proposition \ref{16122011-1119} and Proposition \ref{16122011-1117} (v), $J_1$ is completely continuous. In the sequel, we may assume that  $(w_n)$ converges almost everywhere to some $w_0\in M_p \setminus \{ 0 \}$ 
and that one has $g\in X_p$ such that $|w_n|\leq g$ a.e. on $\Omega$. Further, note that if $w_n (x)\neq 0$, then
$$
\frac{f(x,s_n w_n(x))}{s_{n}^{p-1}} = \frac{f(x, s_n w_n (x))}{(s_n w_n (x))^{p-1}} (w_n(x))^{p-1} \to 
\rho_0 (x) (w_0(x))^{p-1} \mbox{ as } n\to +\infty,
$$
which, by the dominated convergence theorem, implies that $s_{n}^{1-p} H(s_n w_n, t_n)\to \rho_0 N_p(w_0)$ in $X_p^{*}$.
Hence, a passage to the limit in (\ref{16122011-1248}) yields 
$w_0 = J_1 (N_p (w_0)+ \rho_0 N_p (w_0) )$, i.e. $-A_p w_0 + \rho_0 N_p (w_0)=0$. This is a contradiction 
due to Lemma \ref{11122011-0003} and (\ref{16122011-1303}) is proved.\\
\indent Clearly, (\ref{16122011-1248}) allows us to use the homotopy invariance -- Theorem \ref{16122011-1102} (iii)
to see that $\Deg_{M_p} (A_p, F, B_{M_p} (0,R)) = \Deg_{M_p} (A_p, \rho_0 N_p, B_{M_p} (0,R))$, which together with Theorem \ref{02122011-1924} provides the required formula.\\
\indent (ii) The proof is analogical to that for part (i) and it is left to the reader.  \hfill $\square$

\noindent {\bf Proof of Theorem \ref{29092011-1553}}: Let $\delta>0$ and $R>\delta$ be like in Theorem \ref{16122011-1341}.
Then by use of the additivity property -- Theorem \ref{16122011-1102} (ii), we get
\begin{eqnarray*}
\Deg_{M_p} (A_p, F, B_{M_p} (0,R) \setminus \overline{B_{M_p} (0,\delta)}) 
=  \Deg_{M_p} (A_p, F, B_{M_p} (0,R))\! -\! \Deg_{M_p} (A_p, F, B_{M_p} (0,\delta)) \\
= \left\{ 
\begin{array}{cl}
1, & \mbox{ if } \ \rho_0 (x) > \lambda_{1,p} > \rho_\infty (x)  \ \mbox{ for a.e. } x\in \Omega,\\
-1, &  \mbox{ if } \ \rho_0 (x) < \lambda_{1,p} < \rho_\infty (x) \ \mbox{ for a.e. } x\in \Omega.
\end{array}
\right.
\end{eqnarray*}
Hence the existence property of the topological degree gives the existence of $u\in B_{M_p} (0,R) \setminus \overline{B_{M_p} (0,\delta)}$ such that $A_p (u) = F(u)$, which is a required nonzero nonnegative weak solution of (\ref{29102011-1511}). \hfill $\square$

\end{document}